\documentstyle[11pt,twoside]{article}
\newcommand{\BE}{\begin{equation}}
\newcommand{\EE}{\end{equation}}
\newcommand{\BEA}{\begin{eqnarray}}
\newcommand{\EEA}{\end{eqnarray}}

\begin{document}

\hfill LMU-TPW 98/15 

\hfill KA-TP-15-1998 
\vspace{15mm}

{\Huge{\noindent The Relations of Inner and Outer}} 
\vspace{3mm}

{\Huge{\noindent Differential Calculi on Quantum }}
\vspace{2mm}

{\Huge{\noindent Groups}}
\vspace{12mm}

\noindent {\large {\bf Peter Zweydinger}}
\vspace{7mm}

\noindent {\it Institut f{\"u}r Theoretische Physik der Universit{\"a}t (TH) Karlsruhe, Kaiserstr.12 ,
D-76128 Karlsruhe (Germany), Postfach 6980 

\noindent and Max-Planck-Institut f{\"u}r Physik (Werner-Heisenberg-Institut) M{\"u}nchen, F{\"o}hringer 
Ring 6 , D-80805 M{\"u}nchen (Germany)

\noindent and Sektion Physik, Ludwig-Maximilians-Universit{\"a}t M{\"u}nchen, Theresienstr. 37, 
\newline D-80333 M{\"u}nchen (Germany)} 
\vspace{15mm}

\noindent {\bf Abstract} The differential caluli $(\Gamma  , d)$ on 
quantum groups are classified due to the property of the generating 
element $X$ of its differential $d$. There are, on the one hand, 
differential caluli which contain this element $X$ in the basis of 
one- forms that span $\Gamma \;$, called Inner Differential Calculi. On 
the other hand, one has the differential caluli which do not contain the 
generating element $X$ of its differential $d$; thus they are called 
Outer Differential Calculi. We show that this two classes of 
differential caluli, for a given quantum group ${\cal A}$, are related 
by homomorphisms, which map the elements of one class onto elements of 
the other class.
\vspace{15mm}

\noindent {\Large {\bf 1.Introduction}}
\vspace{6mm}

The differential calulus on quantum groups, as a special branch of 
Noncommutative Geometry in the sense of A.Connes \cite{Co1}\cite{Co2}, 
should play a prominent role in the construction of physical models 
based on quantum groups. This stems from the general concept of physics 
that the dynamics (or sometimes the kinematics), i.e. the equations of 
motion, are expressed as (systems of) differential equations. Thus one 
expects in the noncommutative setting the same, or at least a natural 
counterpart of such structures.

The abstract structure of quantum groups is that of non-commutative 
non-cocommuta-
\linebreak 
tive Hopf algebras. The objects of our consideration are 
more precisely the noncommutative generalizations of the 
functionalgebras $Fun(G)$ of some Lie group $G$, which we will call for 
short quantum groups. General references on quantum groups are for 
example the books \cite{KS}\cite{Ks}.

The first example of a differential calculus on quantum groups was given 
by Woronowicz \cite{Wo1}. Shortly afterwards he gave in \cite{Wo2} the 
general theory of differential calculi on quantum groups. Later there 
were, and still are, many contributions from several authors to this 
field. A necessarily personal selection of such contributions is 
Jur$\check{\rm c}$o \cite{Ju}, Carow-Watamura/Schlieker/Watamura/Weich 
\cite{CSWW},  Rosso \cite{Ro}, Aschieri/Castellani \cite{AsCa}, 
Schm{\"u}dgen/Sch{\"u}ler \cite{SS1}\cite{SS2}. What can be 
extracted from all this papers is, as far as bicovariant differential 
calculi are concerned, that there are two classes of differential 
calculi. We will restrict ourself to this situation, because bicovariant 
differential calculi are the most natural generalizations of the 
commutative  differential calculus.

The above mentioned classification of the bicovariant differential 
calculi $(\Gamma  , d)$ is due to the properties of the generating 
element $X$ of the exterior differential $d$. On the one hand we have the 
possibility, that the generating element $X$ does not belong to 
$\Gamma$, thus we will call such differential calculi Outer Differential 
Calculi\footnote{In the sense of an outer automorphism.}. On the other 
hand there are Inner Differential Calculi\footnote{In the sense of an 
inner automorphism.}, which are characterized by the fact that the 
generating element $X$ of its differential $d$ is an element of the set 
of one-forms $\omega^{i}$, which span $\Gamma$ (i.e. $X \in \Gamma$ ).

In order to show how these two classes are related, we will proceed as 
follows. Firstly we will give a review of the construction of  Outer 
Differential Calculi, which is nothing else but the original Woronowicz 
construction \cite{Wo2}. Afterwards we will give a short description of 
the Inner Differential Calculi, as it was initiated independently by 
Jur$\check{\rm c}$o \cite{Ju} and 
Carow-Watamura/Schlieker/Watamura/Weich \cite{CSWW}. This differential 
calculi\footnote{Observe that in this context there was used for the 
first time, as far as I know, the term inner differential calculi, thus 
implying implicitly our splitting.} are classified by 
Schm{\"u}dgen/Sch{\"u}ler \cite{SS1}\cite{SS2} for the q-analogues of the 
standard series of Lie groups. This review is mainly intended to make 
the paper selfcontained, and to fix the notation which will be used in 
the following. All this will be done in the section 2.

In section 3 we show, starting from the splitting of bicovariant 
differential calculi into Outer Differential Calculi and Inner 
Differential Calculi, that using one-dimensional bicovariant 
differential calculi, it is possible to map from one class into the other 
and vice versa. This is done by construction of the maps 
$\Delta^{Out}_{In}$ and $\Delta^{In}_{Out}$ which map the set of 
isomorphy-classes $\{(\Gamma ,\partial)\}$ of Outer Differential 
Calculi into the set of isomorphy-classes $\{(\hat{\Gamma} , d)\}$ of 
Inner Differential Calculi and vice versa. Furthermore we give 
explicitly maps, which associate to any differential calculus in one 
class one in the other class, and show some of their properties.

In section 4 we give, as an application the decomposition of the q-de 
Rham complex of the Inner Differential Calculi in the differential 
bicomplex $(\Gamma^{r,s} , \partial ,\delta)$, which consists of the 
q-de Rham complex of the Outer Differential Calculi and the q-de Rham 
complex of the one-dimensional bicovariant differential calculi whose 
relation was shown in section 3.
\vspace{10mm}

\noindent {\Large {\bf 2.The Bicovariant Differential Calculi on Quantum 
Groups }}
\vspace{6mm}

The bicovariance of differential calculi on quantum groups is the 
substitute of the commutativity of the left and right action of 
differential forms in the case of Lie groups. So they are in a certain 
sense the most natural differential calculi on quantum groups. In 
subsection 2.1 we will give the general properties of bicovariant 
differential calculi, while describing the Outer Differential Calculi, 
which is nothing else but a short version of the original Woronowicz 
approach \cite{Wo2}. In subsection 2.2 we will give a short review of the 
properties which distinguish the Inner Differential Calculi.

Having a look on the respective constructions of Inner and Outer 
Differential Calculi, we can see that there is also the possibility to 
give the name cotangent approach to the Outer Differential Calculi, and 
tangent approach to the  Inner Differential Calculi. 
\vspace{6mm}

\noindent {\large {\bf 2.1 The Outer Differential Calculi on Quantum 
Groups }}
\vspace{4mm}

In this section we shortly review the foundational article \cite{Wo2}, 
where the reader can find the details (and is recommended to do so), 
especially concerning proofs, which wo omit here. Only where we feel it 
absolutely necessity for the understanding will we sketch the proofs. 
Our intention is to be close to the original article, but to stress 
those points which justify the notion of Outer Differential Calculi. In 
order to facilitate a comparison with the original article, we will 
quote in brackets its equation numbers of Theorems, etc. 

We start from the abstract form of the bundle theoretic point of view on 
first order differential forms. (This is the reason why we call this 
calculi cotangent approach differential calculi)

\noindent {\bf Definition 1}[Wo2, Def.1.1]{\sl Let ${\cal A}$ be an 
algebra with unity and ${\Gamma}$ be a bimodule over ${\cal A}$ with
\BE
\partial : {\cal A} \longrightarrow \Gamma  \; .
\EE 
a linear map. We say that  $(\Gamma ,\partial)$ is a first order 
differential calculus over ${\cal A}$ if:
\noindent-1. For any $a,b \in {\cal A}$ the Leibniz-rule holds :
\BE
\partial(ab) = \partial(a)b + a\partial b \; .
\EE
-2. Any element  $\varrho \in \Gamma$ is of the form:
\BE
\varrho = \sum_{k=1}^{K} a_{k}\partial(b_{k})
\EE  
where $a_{k},b_{k} \in {\cal A}$, ($k = 1, \dots , K$ and $K$) is a 
positive integer.}

Two first order differential calculi $(\Gamma ,\partial)$ and 
$(\Gamma' ,\partial')$ are called isomorphic if there exists a 
bimodule-isomorphism ${\rm iso}:\Gamma \rightarrow \Gamma'$ such that 
${\rm iso} \circ \partial a = \partial' a$ for any $a \in {\cal A}$. In 
the following we will identify such isomorphic differential calculi.

For every algebra ${\cal A}$ with the multiplication- map $m : {\cal A} 
\otimes {\cal A} \rightarrow {\cal A}$ there exists a linear subset:
\BE
{\cal A}^{2} = \{ q \in  {\cal A} \otimes {\cal A} | \; mq = 0 \; \}\; .
\EE
of ${\cal A} \otimes {\cal A}$. With the relations 
\BEA
 c\bigg(\sum_{k} a_{k} \otimes b_{k} \bigg) &=& \sum_{k} 
ca_{k} \otimes b_{k} \\
\bigg(\sum_{k} a_{k} \otimes b_{k} \bigg)c &=& \sum_{k} 
a_{k} \otimes b_{k}c
\EEA
for any $ c \in {\cal A}$, ${\cal A}^{2}$ has a 
${\cal A}$-bimodule-structure. Setting, for any  $b \in {\cal A}$ and 
using the unit-element $I$ of ${\cal A}$ , 
$Db = I \otimes b - b \otimes I$ one can show that $D : {\cal A} 
\rightarrow {\cal A}^{2}$ is a linear map and $({\cal A}^{2} , D)$ is a 
first order differential calculus over ${\cal A}$ (If ${\cal A}$ has a 
unity, this is an universal construction.). Its importance stems from 
this universality.

{\bf Proposition 2}[Wo2, Prop.1.1]: {\sl Let ${\cal N}$ be a subbimodul 
of ${\cal A}^{2}$, with $\Gamma = {\cal A}^{2}/{\cal N}$ and  $\pi$ is 
the canonical epimorphism\footnote{An epimorphism is a surjective 
homomorphism.} ${\cal A}^{2} \rightarrow \Gamma$ and $\partial = \pi 
\circ D$. Then $(\Gamma ,\partial)$ is a  first order differential 
calculus over ${\cal A}$. Any first order differential calculus over 
${\cal A}$ can be obtained in this way.}

Up to now these are totally general considerations. For quantum groups we 
have to introduce the following notations: $\phi$ is the 
comultiplication, $\kappa$ is the antipode and $\epsilon$ is the unit. 
The bicovariance of some bimodule $\Gamma$ of a quantum group ${\cal A}$ 
is given by the existence of the two linear maps $\;_{\Gamma}\phi :
\Gamma \rightarrow \Gamma \otimes {\cal A}$ and $ \phi_{\Gamma} :\Gamma 
\rightarrow {\cal A} \otimes \Gamma$. They have to fulfill the conditions
\BEA
\mbox{left covariance:} \hspace{1.5cm}
\begin{array}[t]{l@{\: = \:}l} 
    \phi_{\Gamma}(a\varrho b) & \phi(a)\phi_{\Gamma}(\varrho)\phi(b) \\
    (id \otimes \phi_{\Gamma}) \circ \phi_{\Gamma} & (\phi \otimes id) 
    \circ \phi_{\Gamma} \\ 
    (\epsilon \otimes id)\circ \phi_{\Gamma}(\varrho) & \varrho
\end{array}
\nonumber \\
\mbox{ right covariance:}\hspace{1cm}
\begin{array}[t]{l@{\: = \:}l} 
    \;_{\Gamma}\phi(a\varrho b) & \phi(a)\;_{\Gamma}\phi(\varrho)
\phi(b) \\
    (\;_{\Gamma}\phi \otimes id) \circ \,_{\Gamma}\phi & 
(id \otimes \phi ) \circ  \,_{\Gamma}\phi \\ 
    (id \otimes \epsilon)\circ \,_{\Gamma}\phi(\varrho) & \varrho \; ,
\end{array}
\EEA
for all $\forall \, a,b \, \in {\cal A}$ and $\forall \, \varrho \, \in 
\Gamma$. If these two maps are compatible, that is if they obey the 
relation 
\BE
(id \otimes \,_{\Gamma}\phi ) \circ  \phi_{\Gamma} = (\phi_{\Gamma} 
\otimes id)\circ \,_{\Gamma}\phi 
\EE
the triple $(\Gamma , \;_{\Gamma}\phi \, , \phi_{\Gamma})$ gives a 
bicovariant bimodule.

The application of the bicovariance condition on some differential 
calculus $(\Gamma , \partial)$ over ${\cal A}$ gives us a representation 
of the maps $\;_{\Gamma}\phi \, , \phi_{\Gamma}$: 
\BEA
\phi_{\Gamma}(a\partial b) = &\phi(a)\phi_{\Gamma}(\partial b)& = 
\phi(a)(id \otimes \partial)\phi(b) \\
\,_{\Gamma}\phi(a\partial b) = &\phi(a)\,_{\Gamma}\phi(\partial b)& = 
\phi(a)(\partial \otimes id)\phi(b) \; .
\EEA
An element $\omega \, \in \, \Gamma$ is called left- (respectivly 
right-) invariant if and only if it obeys
\BEA
\phi_{\Gamma}(\omega) &=& (I \otimes \omega) \\
\,_{\Gamma}\phi(\omega) &=& (\omega \otimes I) \; .
\EEA
If both conditions are fulfilled simultaneously $\omega \, \in \, 
\Gamma$ is called biinvariant.

Now one can show that there are sets of left- (respectively right-) 
invariant elements $\,_{{\rm inv}}\Gamma$ (respectively 
$\Gamma_{{\rm inv}}$) of left- (respectivly right-) covariant bimodules 
$\Gamma$ over ${\cal A}$. They form a basis $\omega^{\cal I}$ 
(respectively $\eta^{\cal I}$ ), with $\cal I$ an Indexset, of the 
bimodule $\Gamma$ under consideration. With this basis it can be shown 
that every left- (respectively right-) covariant bimodule $\Gamma$ over 
${\cal A}$ has a unique decomposition. (We give the relations for the 
rightcovariant bimodule in brackets.)
\BEA
\varrho = \sum_{i\in \cal I} a_{i}\omega_{i} & & 
\bigg(\varrho = \sum_{i\in \cal I} a_{i}\eta_{i}\bigg)  \\
\varrho = \sum_{i\in \cal I} \omega_{i} b_{i} & & 
\bigg(\varrho = \sum_{i\in \cal I} \eta_{i} b_{i}\bigg) \; ,
\EEA
which are related by a class of characteristic functionals 
$f_{ij} \, , \, \tilde{f}_{ij} \, \in 
{\cal A}\,'$ (where ${\cal A}\,'$ is the space of linear functionals 
over ${\cal A}$), respectively, according to: 
\BEA
\omega_{i}b  = \sum_{j\in \cal I} (f_{ij}*b)\omega_{j} & & 
\bigg(\eta_{i}b = \sum_{j\in \cal I} (b* \tilde{f}_{ij})\eta_{j}\bigg)  \\
a\omega_{i} = \sum_{j\in \cal I} \omega_{i}((f_{ij} \circ 
\kappa^{-1})*a)  & & 
\bigg(a\eta_{i} = \sum_{j\in \cal I} \eta_{j}(a*(\tilde{f}_{ij} \circ 
\kappa^{-1})\bigg) \; .
\EEA 
where we have used the convolution product $(\xi \ast a)  
\stackrel{def}{=} (id \otimes \xi)\circ \phi(a)$ 
and $(b \ast \xi) \stackrel{def}{=} (\xi \otimes id)\circ \phi(b)$ for 
all $a,b \in {\cal A}$ and $\xi \in {\cal A}\,'$. This characteristic 
functionals have to fulfill the following compatibility relations,
\BEA
f_{ij}(ab) = \sum_{k\in \cal I} f_{ik}(a)f_{kj}(b) & & 
\bigg(\tilde{f}_{ij}(ab) = \sum_{k\in \cal I}  
\tilde{f}_{ik}(a)\tilde{f}_{kj}(b) \bigg) \\
f_{ij}(I) = \delta_{ij} & & \bigg(\tilde{f}_{ij}(I) = \delta_{ij} \bigg)
\EEA
where $a\, , a_{i}\, , b\, , b_{i} \, \in {\cal A}$. 
 
Putting the decomposition and the bicovariance properties together, we 
obtain a relation for the characteristic functionals, which fixes 
them uniquely.

\noindent {\bf Theorem 3}[Wo2, Thm.2.4]: {\sl Let $(\Gamma\, , 
\phi_{\Gamma}\, , \,_{\Gamma}\phi )$ 
be a bicovariant bimodule over ${\cal A}$, and $\{\omega_{i}\, |\, i 
\in {\cal I} \}$ be the basis in the vector space of all left-invariant 
elements of $\Gamma$. Then we have

\noindent 1. For any $i\in \cal I$ 
\BE
\,_{\Gamma}\phi(\omega_{i}) = \sum_{j\in \cal I} \omega_{i} \otimes R_{ji}
\EE
where $R_{ji} \in {\cal A} \,(i,j \in \cal I)$ satisfy the following 
relations
\BEA
\phi(R_{ji}) &=& \sum_{k\in \cal I} R_{jk} \otimes R_{ki} \\
\epsilon(R_{ji}) &=& \delta_{ji} \; . 
\EEA
\noindent 2. There exists a basis $\{ \eta_{i} | i \in \cal I \}$ 
(indexed by the same Indexset $\cal I$) in the vector space of all 
right-invariant elements of $\Gamma$, such that
\BE
\omega_{i} = \sum_{j\in I} \eta_{j}R_{ji}
\EE
is fulfilled for all $i \in \cal I$.

\noindent 3. With this choise of  basis in $\Gamma_{{\rm inv}}$ the 
functionals $f_{ij}$ and $\tilde{f}_{ij}$ with $(i,j \in \cal I)$, 
introduced in eqns. (15)-(18), coincide. 

\noindent 4. For any  $i,j,h \in \cal I$ and any $a \in {\cal A}$ 
we have the relation:
\BE
\sum_{i\in I} R_{ij}(a*f_{ih}) = \sum_{i\in \cal I} (f_{ji}*a)R_{hi}\; .
\EE } 

The construction of the exterior algebra over the bimodule 
$(\Gamma\, , \phi_{\Gamma}\,   \,_{\Gamma}\phi )$ is made as follows. We 
start with a graded algebra of ${\cal A}$-modules 
$T^{n}$ (with ${\cal A} = T^{0}$)
\BE
T =  \sum_{n=0}^{\infty}\!\!^{\bigoplus}\, T^{n} \hspace{2cm} 
\mbox{with} \hspace{1cm} T^{n} = \bigotimes_{i = 0}^{n}\!_{\cal A} \; T 
\; ,
\EE
where the tensor-product of the ${\cal A}$-modules $T$ (respectively 
$\Gamma$) is given by the relation $\omega \otimes_{\cal A} \, 
a \eta = \omega a \, \otimes_{\cal A} \,\eta$ for 
$\forall \; \omega \,  , \eta \; \in T $ (respectively $\Gamma$) and 
$\forall \; a \, \in \,{\cal A}$. Associating the two grade preserving 
linear maps $\phi_{T} : T \rightarrow {\cal A} \otimes T$ and 
$\,_{T}\phi : T \rightarrow T \otimes {\cal A}$, that is:
\BEA
\phi_{T}(T^{n}) &\subset& {\cal A} \otimes T^{n} \\
\,_{T}\phi(T^{n}) &\subset& T^{n} \otimes {\cal A}
\EEA
this makes $(T \, , \phi_{T} \, ,\,_{T}\phi)$ a bicovariant, graded 
algebra. The restriction of $\phi_{T}$ (respectively $\,_{T}\phi$) on 
$T^{n}$ is called $\phi_{T}^{n}$ (respectively $\,_{T}\phi^{n}$). If we 
have $T^{0} = {\cal A}$ and $\phi_{T}^{0} = \,_{T}\phi^{0} = \phi$, then 
$(T^{n} \, , \phi_{T}^{n} \, ,_{T}\phi^{n})$ is a bicovariant 
bimodule over ${\cal A}$.

Let $(\Gamma\, , \phi_{\Gamma}\, , \,_{\Gamma}\phi )$ be a given 
bicovariant bimodule over ${\cal A}$, then we call$(T \, , \phi_{T} \, ,
\,_{T}\phi)$ a graded algebra built over 
$(\Gamma\, , \phi_{\Gamma}\, , \,_{\Gamma}\phi )$ if

\noindent 1. $T^{0} = {\cal A}$ and $\phi_{T}^{0} = \,_{T}\phi^{0} = 
\phi$.

\noindent 2. The bicovariant bimodule $(T^{1} \, , 
\phi_{T}^{1} \, , \,_{T}\phi^{1})$ coincides with $(\Gamma\, , 
\phi_{\Gamma}\, , \,_{\Gamma}\phi )$.

\noindent 3. $T$ is generated by an element of grade $1$, i.e. any 
element $\tau \in T^{n} (n = 2,3, \dots )$ is of the form  
$\tau = \sum_{i} \tau^{i}$, where $\tau^{i}$, for all $i$, is a product 
of $n$ elements of $\Gamma$.

Applying this method it is possible to associate to any bimodule 
$(\Gamma\, , \phi_{\Gamma}\, , \,_{\Gamma}\phi )$ two bicovariant graded 
bimodules $(\Gamma^{\otimes}\, , \phi_{\Gamma}^{\otimes}
\, , \,_{\Gamma}\phi^{\otimes} )$, the tensor-algebra, and 
$(\Gamma^{\wedge}\, , \phi^{\wedge}_{\Gamma}\, , 
\,_{\Gamma}\phi^{\wedge} )$, the exterior algebra. The exterior algebra 
$(\Gamma^{\wedge}\, , \phi^{\wedge}_{\Gamma}\, 
, \,_{\Gamma}\phi^{\wedge} )$ is obtained from the tensor algebra 
$(\Gamma^{\otimes}\, , \phi_{\Gamma}^{\otimes}\, 
, \,_{\Gamma}\phi^{\otimes} )$ due to the factorization of the symmetric 
part (for details see chapter 3 of \cite{Wo2}).

This exterior algebras are characterized as follows:

\noindent {\bf Theorem 4}[Wo2, Thm.3.3]: {\sl Let $(\tilde{\Gamma} 
\, , \tilde{\phi}_{\Gamma}\, , \,_{\Gamma}\tilde{\phi})$ be a 
bicovariant bimodule over ${\cal A}$ and $(\tilde{\Gamma}^{\wedge}
\, , \tilde{\phi}_{\Gamma}^{\wedge}\, , 
\,_{\Gamma}\tilde{\phi}^{\wedge} )$ be the exterior algebra built over
 $(\tilde{\Gamma} \, , \tilde{\phi}_{\Gamma}\, , 
\,_{\Gamma}\tilde{\phi})$, with the
left- and right-invariant sub-bimodule $\Gamma$ of  
$\tilde{\Gamma}$. The maps $\phi_{\Gamma}$ 
(respectively $\,_{\Gamma}\phi $) are restrictions of the maps 
$\tilde{\phi}_{\Gamma}$ ( respectivly $\,_{\Gamma}\tilde{\phi}$) on 
$\Gamma$ and $(\tilde{\Gamma}^{\wedge}\, , \tilde{\phi}_{\Gamma}^{\wedge}
\, , \,_{\Gamma}\tilde{\phi}^{\wedge} )$ on the exterior algebra built 
over $(\Gamma\, , \phi_{\Gamma}\, , \,_{\Gamma}\phi )$. Then there 
exists a gradepreserving multiplicative embedding
\BE
\Gamma^{\wedge} \subset \tilde{\Gamma}^{\wedge} \; ,
\EE
coinciding with $id$ on elements of grade $0$ (i.e. on ${\cal A}$) 
and with the inclusion $\Gamma \subset \tilde{\Gamma}$ on elements of 
grad $1$.Moreover this inclusion intertwines the left and right action 
of ${\cal A}$ on $\Gamma^{\wedge}$ and $\tilde{\Gamma}^{\wedge}$, 
respectively. }

The use of this construction on first order differential calculi gives 
then higher order differential calculi , which are described in the 
following theorem.

\noindent{\bf Theorem 5}[Wo2, Thm.4.1]: {\sl Let $G = (A,u)$ be 
compact quantum group, and $(\Gamma \, , \partial)$ be a first 
order differential calculus on $G$; let $\phi_{\Gamma}$ and 
$\,_{\Gamma}\phi$ be the left- respective right- actions of $G$ on
$\Gamma$ given in eqns.(7), and let 
$(\Gamma^{\wedge}\, , \phi^{\wedge}_{\Gamma}\, , 
\,_{\Gamma}\phi^{\wedge} )$ be the exterior algebra built over 
$(\Gamma\, , \phi_{\Gamma}\, , \,_{\Gamma}\phi )$. Then there exists one 
and only one lineare map
\BE
\partial :\Gamma^{\wedge} \longrightarrow \Gamma^{\wedge} \, ,  
\EE 
such that

\noindent 1. $\partial$ increases the grade by one.
 
\noindent 2. On elements of grade 0 , $\partial$ coincides with the 
original derivation from eqn.(1).

\noindent 3. $\partial$ is a graded derivative:
\BE
 \partial(\Theta \wedge \Theta') = \partial(\Theta) \wedge \Theta' 
+ \, (-)^{k} \Theta \wedge \partial \Theta'
\EE
for any $\Theta \in \Gamma^{\wedge k}$ and $ \Theta'\in 
\Gamma^{\wedge}$ with $k \in {\rm {\bf N}}_{0}$.

\noindent 4. $\partial$ is a coboundary 
\BE
\partial(\partial \Theta ) = 0
\EE
for any $\Theta \in \Gamma^{\wedge}$.

\noindent 5. $\partial$ is bicovariant:
\BEA
\phi^{\wedge}_{\Gamma}(\partial \Theta) &=& (id \otimes \partial)
\phi^{\wedge}_{\Gamma}(\Theta) \\
\,_{\Gamma}\phi^{\wedge}(\partial \Theta) &=& (\partial \otimes id)
\,_{\Gamma}\phi^{\wedge}(\Theta)
\EEA
for any $\Theta \in \Gamma^{\wedge}$.}

\noindent{\bf Proof} (sketch): The applicatin of the extended bimodule 
construction given in chapter 4 of \cite{Wo1}, proofs the theorem.

We associate to the bimodule $\Gamma$ a free left ${\cal A}$-module 
$\Gamma_{0}$,which is generated by one element $X$:
\BE
\tilde{\Gamma} = \Gamma_{0} \bigoplus \Gamma = {\cal A}X \bigoplus \Gamma
\EE
thus any element $\tilde{\xi} \in \tilde{\Gamma}$ has the form 
\BE
\tilde{\xi} = cX + \xi \; ,
\EE
uniquely determined for any  $c \in {\cal A}$ and $\xi \in \Gamma$.

The right-multiplication of elements of ${\cal A}$ on  $\tilde{\Gamma}$ 
is defined by the relation 
\BE
\tilde{\xi}a \stackrel{def}{=} caX + c(\partial a + \xi a) \hspace{2cm}
\forall \, a \in {\cal A} \; .
\EE
From this one can check that $\tilde{\Gamma}$ is a bimodule over 
${\cal A}$ and furthermore, with an appropriate choice (namely 
$(c = I \, , \, \xi = 0)$), one obtains the defining relation of the 
exterior derivative $\partial$:
\BE
\partial a = Xa -aX
\EE 
The fact that $X \notin \Gamma$ is the justification to call such 
calculi {\bf Outer Differential Calculi} (in the sense of an outer 
automorphism).

The choice 
\BEA
\tilde{\phi}_{\Gamma}(\tilde{\xi}) &=& {\phi}(c)(X \otimes I) 
+ {\phi}_{\Gamma}(\xi) \, ,\\
\,_{\Gamma}\tilde{\phi}(\tilde{\xi})) &=& {\phi}(c)(I \otimes X) 
+\; \,_{\Gamma}{\phi}(\xi)
\EEA
of the extended left- (respectively right-) action of $G$ makes 
$(\tilde{\Gamma} \, , \tilde{\phi}_{\Gamma}\, 
, \,_{\Gamma}\tilde{\phi})$ a bicovariant bimodule over ${\cal A}$, 
which contains $(\Gamma\, , \phi_{\Gamma}\, , \,_{\Gamma}\phi )$ as a 
invariant sub-bimodule. The extension of 
$(\tilde{\Gamma} \, , \tilde{\phi}_{\Gamma}\, , \,_{\Gamma}\tilde{\phi})$ 
to the exterior bicovariant bimodule $(\tilde{\Gamma}^{\wedge}\, , 
\tilde{\phi}_{\Gamma}^{\wedge} \, , \,_{\Gamma}\tilde{\phi}^{\wedge})$ 
makes $\Gamma^{\wedge}$ according to theorem 4 an exterior bicovariant 
differential bimodule.

From eqns.(37) and (38) one obtains
\BE
 X\wedge X = 0 \; .
\EE
This allows us to put, as a generalization of eqn. (36) to $\forall 
\Theta \in {\tilde{\Gamma}}$,
\BE
\partial \Theta = \left[ X \, , \,\Theta \right]_{grad} 
\EE 
where the graded commutator $\left[  \;\, , \; \right]_{grad} $ is given 
by
\BE
\left[ X \, , \,\Theta \right]_{grad} \stackrel{{\rm def}}{=}
 \left\{\begin{array}{ll}
       X \wedge \Theta - \Theta \wedge X & \mbox{if grad }\Theta 
       \mbox{ is even} \\
       X \wedge \Theta + \Theta \wedge X & \mbox{if grad }\Theta 
       \mbox{ is odd} \; .
       \end{array}
 \right.
\EE
Thus we have the following mapping property
\BE
\left[  \;\,  , \; \right]_{grad} : \tilde{\Gamma}_{0} \times 
\tilde{\Gamma}^{\wedge} \longrightarrow \Gamma^{\wedge} \; .
\EE
The coboundary condition 
\BE
\partial^{2} = 0 
\EE
is then implied by this construction.
\hfill \rule{3mm}{3mm}
\vspace{3mm}

Now one constructs the vector space of vectorfields $T \subset 
{\cal A}'$ on  ${\cal A}$. It is characterized by the right ideal $\Re$  
of  ${\cal A}$, which is defined as the subset of the kernel of the 
counit $\ker(\epsilon)$ that gives the subbimodule ${\cal N}$ from 
Prop.2 as follows
\BE
T = \left\{\chi \in {\cal A}' \, | \, \chi(I) = 0 \,\; \bigwedge \; 
\chi(a) = 0 \;\; \mbox{for} \;\; \forall \, a \in \Re \right\} \, .
\EE 
Then we have
 
\noindent{\bf Theorem 6}[Wo2, Thm.5.1]: {\sl There exists a unique 
bilinear form:
\BE
 \left< \; , \; \right> : \Gamma \times T \rightarrow {\bf C} 
\hspace{2cm} (\varrho , \chi) \in \Gamma \times T : \, 
 \left< \varrho , \chi \right> \in {\bf C} \; , 
\EE
such that for any $a \in {\cal A}\, ,\varrho \in \Gamma$ and 
$\chi \in T$ we have
\BEA
 \left< a\varrho , \chi \right> &=&  \epsilon(a) \left< \varrho , 
       \chi \right> \\
 \left< \partial a , \chi \right> &=&  \chi(a) \; .
\EEA
Moreover, denoting by $\,_{inv}\Gamma$ the set of all leftinvariant 
elements of $\Gamma$, we have

\noindent 1. For any $\omega \in \,_{inv}\Gamma$, 
\BE
\left( \begin{array}{c}  
         \left< \omega , \chi \right> = 0 \\
         \mbox{for } \forall \, \chi \in T
       \end{array}  
\right) 
\Longrightarrow \left( \omega = 0 \right) \; .
\EE
\noindent 2. For any $\chi \in T$, 
\BE
\left( \begin{array}{c}  
\left< \omega , \chi \right> = 0 \\
        \mbox{for } \forall \, \omega \in \,_{inv}\Gamma 
         \end{array}  
\right) 
\Longrightarrow \left( \chi = 0 \right) \; .
\EE  }
This makes it natural that $T$ and $\,_{inv}\Gamma$ are a dual pair of 
vector spaces with respect to $\left< \; , \; \right>$. Due to this 
bilinear form it is possible to obtain a dual basis $\chi_{i}$ of vector 
fields, which one obtains from the basis $\omega^{i}$ of $\,_{inv}\Gamma$.
\BE
\left< \omega^{i} \; , \; \chi_{j} \right> = \delta^{i}_{\, j} 
\hspace{3cm}
(i,j \in \cal I, \mbox{with } \cal I \mbox{ is an index- set}) \; .
\EE
From eqn.(47) one obtains a family $a_{i} (i \in \cal I)$ of elements of 
$\ker(\epsilon)$ which fulfill 
the relation 
\BE
\chi_{i}(a_{j}) = \delta_{ij} \; .
\EE
They give us a connection of the basis $\omega^{i}$ and the exterior 
differential $\partial$.

\noindent {\bf Theorem 7}[Wo2, Thm.5.2]: 
{\sl 1. For any $a \in {\cal A}$
\BE
\partial a = \sum_{i \in \cal I} (\chi_{i} * a) \omega^{i} \, .
\EE
\noindent 2. Let $\{ f_{ij} \, | \, i,j \in I \}$ be a family of linear 
functionals on ${\cal A}$, $(f_{ij} \in {\cal A}')$, introduced in 
eqns.(15)-(18). Then we have, for any $i \in I$ and $a,b \in {\cal A}$ 
\BE
\chi_{i}(ab) =\sum_{j \in \cal I} \chi_{j}(a)f_{ji}(b) + 
\epsilon(a)\chi_{i}(b) \, .
\EE
In particular we have for any $i,j \in \cal I$ and $b \in {\cal A}$
\begin{center}
$\chi_{i}(a_{j}b) = f_{ji}(b)$~.
\end{center}
\noindent 3. Let $\{ R_{ij} \, | \, i,j \in \cal I \}$ the family of 
elements of  ${\cal A}$ introduced in Theorem 3. Then
\BE
R_{ij} = (\chi_{i} \otimes id)({\rm ad}(a_{j}))  \, .
\EE  }

Here we have used the definition $ {\rm ad}(a) = \sum_{k} b_{k} \otimes 
\kappa(a_{k})c_{k}$, with $a_{k},b_{k},c_{k} \in {\cal A}$ such that 
$(id \otimes \phi)\circ \phi =  \sum_{k} a_{k} \otimes b_{k} \otimes 
c_{k}$. 

Due to the structure of the exterior algebra $\Gamma^{\wedge}$ it is 
possible to construct a higher order differential calculus. The 
q-Lie-bracket is then given by
\BE
\left[  \chi' \, , \, \chi'' \right] \stackrel{{\rm def}}{=} 
\chi' * \chi'' - \sum_{s} \chi_{s}'' * \chi_{s}' \;\; ,
\EE  
with $\chi_{s}' , \chi_{s}''\, (s = 1, \dots , S)$ and $S \in 
\, {\bf N}$ are elements of $T$, such that\footnote{$\sigma^{t}$ is the 
twistmap in the space $T \otimes T$ dual to the twist-map 
$\sigma_{\Gamma}$ on the space $\Gamma \otimes \Gamma$.} 
\BE
\sigma^{t}(\chi' \otimes \chi'') = \sum_{s} \chi_{s}'' \otimes 
\chi_{s}' \; .
\EE
Then we have 

\noindent {\bf Theorem 8}[Wo2, Thm.5.3]:{\sl 1. For $\forall \chi' , 
\chi'' \in T$
\BE
\left[  \chi' \, , \, \chi'' \right] \in T \, .
\EE
2. If $\chi_{s}' , \chi_{s}'' (s = 1, \dots , S)$ 
are elements of $T$ such that $\sigma^{t}( \sum_{s} \chi_{s}' 
\otimes \chi_{s}'') =
\sum_{s} \chi_{s}' \otimes \chi_{s}''$ then
\BE
\sum_{s} \left[ \chi_{s}' \, , \, \chi_{s}'' \right] = 0 \, .
\EE
3.  For $\forall \chi , \chi' , \chi'' \in T$ 
\BE
\left[  \chi \, ,\left[  \chi' \, , \, \chi'' \right] \right] = 
\left[ \left[  \chi \, , \, \chi' \right], \, \chi'' \right] - 
\sum_{s} \left[ \left[ \chi \, , \,\chi_{s}''  \right]
 , \, \chi_{s}' \right]
\EE
where $\chi_{s}' , \chi_{s}'' (s = 1, \dots , S)$ are elements of 
$T$ introduced in eqn.(56).  }

We now have that eqn.(58) is the q-deformed antisymmetry of the 
commutator, while eqn.(59) gives the q-analog of the Jacobi-identity.

In terms of the basis elements $\omega^{i}$ of $\,_{inv}\Gamma$ this 
reads
\BE
\sigma_{\Gamma}(\omega_{i} \otimes \omega_{j}) = 
\sum_{k.l \in \cal I} \lambda_{ij}^{\;\;\; kl}
(\omega_{k} \otimes \omega_{l}) 
\EE 
with $\lambda_{ij}^{\;\;\; kl} \in {\bf C} \, (i,j,k,l \in I)$,  which 
implies on the dual space $T$ the corresponding relations 
\BE
\sigma^{t}(\chi_{k} \otimes \chi_{l}) = 
\sum_{i,j \in \cal I} \lambda_{ij}^{\;\;\; kl}(\chi_{i} \otimes 
\chi_{j}) \; .
\EE
Thus we obtain for the q-Lie-bracket in terms of the dual basis $\chi_{i}$
\BE 
\left[\chi_{k}  \, , \, \chi_{l} \right] = \chi_{k} * \chi_{l} - 
\sum_{i,j \in \cal I} \lambda^{ij}_{\;\;\; kl}(\chi_{i} * \chi_{j}) \; ,
\EE
and the corresponding form of the q-Jacobi-identity
\BE
\left[\chi_{i}  \, , \left[\chi_{j}  \, , \, \chi_{k} \right] \right] = 
\left[\left[\chi_{i}  \, ,\chi_{j} \right] \, , \, \chi_{k} \right] - 
\sum_{l,m \in \cal I} \lambda^{lm}_{\;\;\;\, jk} 
\left[\left[\chi_{i}  \, ,\chi_{l} \right] \, , \, \chi_{m} \right] \, .
\EE 

The same construction works as well for right invariant differential 
calculi.
\vspace{6mm}

\noindent {\large {\bf 2.2 The Inner Differential Calculi on Quantum 
Groups }}
\vspace{4mm}

The Inner Differential Calculi on quantum groups are based on the 
FRT-construction \cite{FRT} of quantum groups, which uses the duality of 
the algebra of functions (quantum groups) ${\rm A}_{{\rm R}}$ and the 
Universal Envelopping Algebra (QUEA) ${\rm U}_{{\rm R}}$. The regular 
functionals $(L^{\pm})^{a}_{\; b}$ form a basis of the QUEA 
${\rm U}_{{\rm R}}$. Relying on this structure \cite{Ju},\cite{CSWW} 
constructed bicovariant differential calculi on the quantum groups, 
generalizing the $A_{n}-,B_{n}-,C_{n}-,$ $D_{n}-$series of Lie groups. 
We will in our presentation follow the paper \cite{AsCa}, where the 
reader can find the details.

The regular functionals $(L^{\pm})^{a}_{\; b}$ have the characteristic 
contractions with the generators $t^{c}_{\; d}$ of the function-algebra 
$A_{R}$, 
which gives the fundamental R-matrix
\BE
(L^{\pm})^{a}_{\; b}(t^{c}_{\; d} ) = 
R^{\pm ac}_{\;\;\;\;\;\;\; bd} \, .
\EE 
This fixes the QUEA ${\rm U}_{{\rm R}}$. 
(The specific form of the R-matrices can be found in \cite{AsCa}).

The coalgebra-structure of ${\rm U}_{{\rm R}}$, in this representation, 
is given by the following relations (in the following we will use the 
Einstein-summation-convention),
\BEA
\phi_{d}((L^{\pm})^{a}_{\; b}) &=& (L^{\pm})^{a}_{\; g} \otimes 
(L^{\pm})^{g}_{\; b} \\
\epsilon_{d}((L^{\pm})^{a}_{\; b}) &=& \delta^{a}_{\; b} \\
\kappa_{d}((L^{\pm})^{a}_{\; b}) &=& (L^{\pm})^{a}_{\; b} \circ \kappa_{c}
\EEA
where $\phi_{d}$ is the comultiplcation in ${\rm U}_{{\rm R}}$, 
$\epsilon_{d}$ is the counit in ${\rm U}_{{\rm R}}$ and $\kappa_{d}$ is 
the antipode in ${\rm U}_{{\rm R}}$ (the subskripts d (respectively c) 
refer to the dual ${\rm U}_{{\rm R}}$-Hopf-algebra (respectively the 
original function algebra 

$A_{R}$-Hopf-algebra)) structures. This Hopf-algebra properties fixes the 
relations among the generators of the differental calculus.

Due to the consistency condition of the exterior differential $d$ and the 
bimodule-structure, one obtains relations among the vectorfields 
$\chi_{i}$ and the characteristic functionals $f^{i}_{\: j}$ of the 
bimodule (compare this with eqns.(15)-(18) and (44),(50)). These 
bicovariance conditions are
\BEA
\left[\chi_{i}  \, , \, \chi_{j} \right] &\equiv & \chi_{i} \chi_{j} - 
 \Lambda^{kl}_{\;\;\; ij}(\chi_{k} \chi_{l}) = 
{\rm C}_{ij}^{\:\;\; k}\chi_{k} \\
\Lambda^{nm}_{\;\;\; ij} f^{i}_{\: p} f^{j}_{\; q} &=& 
f^{n}_{\: i} f^{m}_{\;\; j} \Lambda^{ij}_{\;\;\; pq} \\
{\rm C}_{mn}^{\:\;\;\; i} f^{m}_{\;\; j} f^{n}_{\;\; k} + 
f^{i}_{\;\; j}\chi_{k} &=& 
\Lambda^{pq}_{\;\;\;\; jk} \chi_{p}f^{i}_{\;\; q} + 
{\rm C}_{jk}^{\:\;\; l} f^{i}_{\;\; l} \\
\chi_{k} f^{n}_{\;\; l} &=& 
\Lambda^{ij}_{\;\;\; kl}f^{n}_{\;\; i} \chi_{j}  \, .
\EEA
Here are the ${\rm C}_{ij}^{\:\;\; k}$ the q-structure constants (which 
are implicit in eqn.(57)) and the $\Lambda^{kl}_{\;\;\; ij}$ are the 
q-commutation functions of the vectorfields (like in eqn.(62)),
(For convenience we have changed the notation of the 
$\lambda^{kl}_{\;\;\; ij}$ from eqn.(62) to $\Lambda^{kl}_{\;\;\; ij}$). 
The Hopf-algebra costructure of ${\rm U}_{{\rm R}}$ expressed in the 
vectorfields $\chi_{i}$ and the characteristic functionals 
$f^{i}_{\: j}$ gets the following form
\BEA
\phi_{d}(\chi_{i}) &=& \chi_{j} \otimes f^{j}_{\;\; i} + 
\epsilon_{c} \otimes \chi_{i} \hspace{2.5cm} 
(\mbox{observe:}\; \epsilon_{c} = I_{d}) \hspace{1cm} \\
\epsilon_{d}(\chi_{i}) &=& 0 \\
\kappa_{d}(\chi_{i}) &=& -\chi_{j}\kappa_{d}(f^{j}_{\;\; i}) \\
\phi_{d}(f^{i}_{\;\; j}) &=& f^{i}_{\;\; k} \otimes f^{k}_{\;\; j} \\
\epsilon_{d}(f^{i}_{\;\; j}) &=& \delta^{i}_{\;\; j} \\
\kappa_{d}(f^{k}_{\;\; j}) f^{j}_{\;\; i} &=& 
\delta^{\;\; k}_{i}\epsilon_{c} = f^{k}_{\;\; j}\kappa_{d}
( f^{j}_{\;\; i}) \, .
\EEA

The functionals $\chi_{i} , f^{i}_{\;\; j} , 
{\rm C}_{ij}^{\:\;\; k} , \Lambda^{ij}_{\;\;\; kl}$, which fulfill the 
bicovariance conditions, can be rephrased in terms of the regular 
functionals  $(L^{\pm})^{a}_{\; b}$ and the R-matrices 
$R^{\pm ac}_{\;\;\;\;\;\;\; bd}$. This can be done with the help of the 
eqns.(65)-(67) and (72)-(77). There one has used the following convention 
$\{\,_{a}^{\;\; b} \hat{=} \,^{i} , \,_{\;\; a}^{b} \hat{=} \,_{i}\}$ 
for the replacement of the ``vector indices'' by ``spinor indices'', i.e. 
double indices which reflect the matrix structur of the $L^{\pm}$- and 
$R^{\pm}$-functionals. Now, using the normalization constant $\lambda$ 
(the simplest choice of which is $q - q^{-1}$), one obtains 
\BEA
\chi_{i} &\equiv& \chi_{\;\; c_{2}}^{c_{1}} = \frac{1}{\lambda}
\left\{ \kappa_{d}((L^{+})_{\;\; b}^{c_{1}})(L^{-})_{\;\; c_{2}}^{b} - 
\delta_{\;\; c_{2}}^{c_{1}}\epsilon_{c} \right\} \\
\Lambda^{ij}_{\;\;\; kl}  &\equiv& 
\Lambda_{a_{1}\;\; d_{1}}^{\;\; a_{2}\;\; d_{2}}|
_{\;\; c_{2}\;\; b_{2}}^{c_{1}\;\; b_{1}} = 
d^{f_{2}}d^{-1}_{\;\;\; c_{2}} R^{f_{2}b_{1}}_{\;\;\;\;\; c_{2}g_{1}} 
(R^{-1})^{c_{1}g_{1}}_{\;\;\;\;\; e_{1}a_{1}}
(R^{-1})^{a_{2}e_{1}}_{\;\;\;\;\; g_{2}d_{1}}
R^{g_{2}d_{2}}_{\;\;\;\;\; b_{2}f_{2}} \hspace{6mm}  \\ 
{\rm C}_{ij}^{\:\;\; k} &\equiv& 
{\rm C}_{\;\; a_{2}\;\; b_{2}}^{a_{1}\;\; b_{1}}|_{c_{1}}^{\;\;\; c_{2}} =
-\frac{1}{\lambda}\left\{ 
\Lambda_{a\;\; d_{1}}^{\;\; a\;\; d_{2}}|
_{\;\; c_{2}\;\; b_{2}}^{c_{1}\;\; b_{1}} - \delta_{\;\; b_{2}}^{b_{1}} 
\delta_{\;\; d_{1}}^{c_{1}}\delta_{c_{2}}^{\;\; d_{2}} \right\} \\
f^{i}_{\;\; j} &\equiv& f_{a_{1}\;\;\;\;\;\; b_{2}}^{\;\;\; a_{2}b_{1}} = 
\kappa_{d}((L^{+})^{b_{1}}_{\;\; a_{1}})(L^{-})^{a_{2}}_{\;\; b_{2}} \, ,
\EEA
with the notation
\begin{center}
$ d^{a} \stackrel{{\rm def}}{=} \left\{
  \begin{array}{l}
        q^{2a-1} \hspace{1cm} \mbox{for the} \; A_{n-1}-
        \mbox{series} \\
        D_{\;\; a}^{a} \hspace{1.2cm} \mbox{for the} \; B_{n}, C_{n}, 
        D_{n}-\mbox{series, with} \;D_{\;\; a}^{a} = C^{ae}C_{ae} 
  \end{array}
\right.$
\end{center} 
the details of which can be found in \cite{AsCa}. So we have the tools 
at hand to construct a bicovariant differential calculus from the 
$L^{\pm}$-representation of ${\rm U}_{{\rm R}}$, which is the reason 
to call it a tangent approach.

The exterior differential calculus\footnote{In the sense of exterior 
algebra}, dual to this construction, is obtained by the application of 
the bilinear form $\left< \; , \; \right>$ given in theorem 6. We define 
the $A_{R}$-bimodule $\hat{\Gamma}$ by the basis $\omega_{a}^{\;\; b}$ of
differential-one-forms ($a,b = 1, \dots ,N$), where $N$ is the dimension 
of the fundamental representation of $A_{R}$. This basis is assumed to 
be left invariant\footnote{Analogously one can give right invariant 
calculi}, that is $\,_{inv}{\hat{\Gamma}}$ has dimension $N^{2}$ and we 
have (like in eqn.(11))
\BE
\phi_{\hat{\Gamma}}(\omega_{a}^{\;\; b}) = I \otimes \omega_{a}^{\;\; b}
\hspace{2cm} ,  \hspace{2cm} \omega_{a}^{\;\; b} \in 
\,_{inv}{\hat{\Gamma}} \, ,
\EE
which fixes $\phi_{\hat{\Gamma}}$ on the whole of ${\hat{\Gamma}}$. With 
the characteristic functionals (of the (${\rm A}_{R}$ - 
${\hat{\Gamma}}$)-commutation) $f^{i}_{\; j}$, as given in eqn.(81), we 
have 
\BE
\omega_{a_{1}}^{\;\; a_{2}}a = 
(f_{a_{1}\;\;\;\;\;\; b_{2}}^{\;\;\; a_{2}b_{1}} * a)
\omega_{b_{1}}^{\;\; b_{2}} \, .
\EE

To obtain, in this formalism, the exterior algebra ${\hat{\Gamma}}$, we 
need the exterior multiplication of the elements $\omega_{a}^{\;\; b}$. 
This is, in the case of the $B_{n}-,C_{n}-,D_{n}-$ series, for example, 
given by 
\BE
\omega_{i} \wedge \omega_{j} = -Z^{kl}_{\;\;\; ij}\omega_{k} \wedge 
\omega_{l} \, ,
\EE
where we have used the identity 
\BE
Z^{ij}_{\;\;\; kl} \equiv \frac{1}{q^{2} - q^{-2}} \left[
\Lambda^{ij}_{\;\;\; kl} - (\Lambda^{-1})^{ij}_{\;\;\; kl}
\right] \; .
\EE
Here we have, for notational simplicity, returned to the 
``vector indices''. This reflects the generalized commutation properties 
of the quantum group.

The exterior differential $d : \hat{\Gamma}^{\wedge k} \rightarrow 
\hat{\Gamma}^{\wedge k + 1}$ is constructed similar to theorem 5 and the 
eqns.(34),(35) as well as (40),(41), due to the biinvariant element 
$\hat{X} = \sum_{a} \omega_{a}^{\;\; a}$. But observe, in contrary to the 
considerations given there, we have 
\BE
\hat{X} = \sum_{a} \omega_{a}^{\;\; a} \in \hat{\Gamma} \; , 
\hspace{3cm} \left< \hat{X} \right> = \hat{\Gamma}_{0} \, . 
\EE
Thus we have {\bf not} an extended bimodule construction. Nevertheless, 
we have $\forall \, \Theta \in \hat{\Gamma}^{\wedge k}$, analogously to 
eqns.(40),(41)
\BE
d \,\Theta \stackrel{{\rm def}}{=} \frac{1}{\lambda} 
\left[\hat{X} \, , \, \Theta \right]_{grad} = 
\frac{1}{\lambda}\left( \hat{X} \wedge \Theta - 
(-)^{k} \Theta \wedge \hat{X} \right) \; ,
\EE
with the graded commutator
\begin{center}
$\left[ \; , \; \right]_{grad} :\hat{\Gamma}_{0} \times 
\hat{\Gamma}^{\wedge} \longrightarrow \hat{\Gamma}^{\wedge}$~, 
\end{center}
and the normalization constant $\lambda$ as in eqn.(78) 
(respectively (81)). Especially in the case $\hat{\Gamma}^{\wedge 0} = 
{\rm A}_{R}$ we obtain 
\BE
da = \frac{1}{\lambda} \left[\hat{X} \, , \, a \right] \hspace{2cm}
 \forall a \; \in {\rm A}_{R} \, ,
\EE
which obeys the Leibniz rule
\BE
d(ab) = (da)b + adb \; .
\EE

Due to the above construction we have a differential calculus 
$(\hat{\Gamma} , d )$ over ${\rm A}_{R}$, which is different from the 
differential calculi obtained in theorem 5, in so far as the generating 
element $\hat{X}$ of the exterior differential $d$ is contained in 
$\hat{\Gamma}$. This is the reason to call these differential calculi 
{\bf Inner Differential Calculi} (in the sense of an inner automorphism).
\vspace{15mm}

\noindent {\Large {\bf 3. The Reconstruction Theorem}}
\vspace{6mm}

The reconstruction theorem gives the relation between the Outer- and 
Inner Differential Calculi on quantum groups, described in the previous 
section. In order to find such a relation we start by a comparison of 
these two classes of differential calculi on the level of their 
respective structure of differential forms. That is, we consider 
their bicovariant bimodules and the associated exterior differentials,
respectively.

As we have seen in section 2.1, following the original Woronowicz-
construction of the differential calculi $(\Gamma , \partial)$, the 
Outer Differential Calculi on quantum groups are characterized by the 
bicovariant bimodule $\Gamma$ and the exterior differential $\partial$, 
which is explicit in theorem 5. There one finds that, in order to define 
the exterior differential $\partial$ one has to apply the extended 
bimodule construction, i.e. one has to complete the bimodule $\Gamma$ 
with a bimodule $\Gamma_{0}$, which is generated by the element $X$, to 
obtain the extended bimodule $\tilde{\Gamma}$. As a consequence of 
eqns.(33)-(36) the differential $\partial$ is implied by the element 
$X$, which is, due to this fact, called the generating element (of 
$\partial$) (see also the eqns.(40),(41) for the higher differential 
calculi).

The first order differential calculi $(\hat{\Gamma} , d )$ of the class 
of Inner Differential Calculi, on the other hand, as is shown in 
section 2.2, are given due to the bicovariant bimodule $\hat{\Gamma}$ and 
the exterior differential $d$. In contrast to the construction of the 
Outer Differential Calculi, the exterior differential $d$ is 
characterized by the canonical element 
$\hat{X} = \sum_{a} \omega_{a}^{\;\; a}$, which is an element of the 
bicovariant bimodule $\hat{\Gamma}$. (Compare with the eqns.(86),(88) 
and for the higher differential calculi eqn.(87)).

The application of these two constructions of differential calculi on the 
same quantum group ${\cal A}$ leads us, as was shortly sketched above, 
to two different results, namely to $(\hat{\Gamma} , d )$ in the case of 
Inner Differential Calculi and to $(\Gamma , \partial)$ in the case of 
Outer Differential Calculi. Now comparing the vector space dimension of 
the respective bimodules of differential-one-forms of the both classes, 
under the assumption that the subspaces orthogonal to the generating 
element of the respective differential are isomorphic, we find them 
different. This is obviously due to to fact that the generating element 
$\hat{X}$ of the differential $d$ is contained in $\hat{\Gamma}$, 
contrary to $\Gamma$ which does not contain the generating element $X$ 
of its differential $\partial$. On the other hand, the bimodule $\Gamma$ 
is imbedded in the extended bimodule $\tilde{\Gamma}$ 
\begin{center}
$\tilde{\Gamma} = \Gamma \oplus \Gamma_{0} 
\supset \Gamma $~,
\end{center}
that contains the generating element $X$ of $\partial$ (compare 
chapter 4 in \cite{Wo1}). This observation makes it natural to try for an 
identification of the bimodules $\tilde{\Gamma}$ and $\hat{\Gamma}$. 
The deeper reason for such an identification is due to the fact that it 
is always possible to find a $\Gamma$ which is isomorphic to 
$\hat{\Gamma}/\hat{\Gamma}_{0}$. This is true because of the universal 
character of Prop.2. Furthermore the vector space dimension of 
$\tilde{\Gamma} = \Gamma \oplus \Gamma_{0}$ and $\hat{\Gamma}$ coincide. 
Starting on the contrary from any $\Gamma$ we have analogously a 
subbimodule $\hat{\Gamma}/\hat{\Gamma}_{0}$ isomorphic to $\Gamma$. All 
this is complemented by the fact that the generators of 
$\hat{\Gamma}_{0}$ and $\Gamma_{0}$ obey the same relations, and 
furthermore, as will be shown in Lemma 10, there is a possibility to make 
a completion of $\Gamma_{0}$ to a one-dimensional bicovariant 
differential calculus $(\Gamma_{0} , \delta)$.

Having considered the bimodule part of the definition of the two classes 
of differential calculi, we have now to turn to the differentials. The 
exterior differentials $\partial$, respectively $d$, show also different 
properties. While on the one hand the Image ${\rm Im}(\partial)$ of the 
differential  $\partial$, looked upon as a map, does not contain terms 
proportional to its generator $X$ (see eqn.(42))
\BE
\Gamma_{0} \not\subset {\rm Im}(\partial) = \Gamma \; ,
\EE
we have on the other hand that the image ${\rm Im}(d)$ of the 
differential  $d$, generated by $\hat{X}$, contains its generator (see 
eqn.(86)), thus we have
\BEA
{\rm Im}(\partial) &\subseteq& {\rm Im}(d) \nonumber\\
  &\wedge & \\
{\rm Im}(\partial) &\neq& {\rm Im}(d)\nonumber\; .
\EEA  

In order to investigate the relations of the two classes of differential 
calculi it is useful to formalize the decomposition of the bimodules 
$\hat{\Gamma}$ and $\tilde{\Gamma}$ in the parts that generate the 
respective differentials\footnote{In the case of $\tilde{\Gamma}$ the 
generated differential is that of $\Gamma$, i.e. $\partial$.}, that is 
$\Gamma_{0}$ respectively $\hat{\Gamma}_{0}$, and the subspaces which 
complement them.

{\bf Definition 9: }{\sl Let $J$ be the projection operator on the 
bimodules $\, \hat{\Gamma}$, respectively $\tilde{\Gamma}$
\BE
J: \left\{ \begin{array}{l}
           \hat{\Gamma} \longrightarrow \hat{\Gamma}_{0} \\
           \tilde{\Gamma} \longrightarrow \Gamma_{0} \; ,
           \end{array}
   \right.
\EE
with the projection property $J^{2} = J$, and $J^{\bot}$ the 
complementary operator:
\BE
J^{\bot} = \left\{ \begin{array}{l}
                   (id_{\hat{\Gamma}} - J):\hat{\Gamma} \longrightarrow 
                   \hat{\Gamma} / \hat{\Gamma}_{0} \\
                   (id_{\tilde{\Gamma}} - J):\tilde{\Gamma} 
                   \longrightarrow \tilde{\Gamma} /\Gamma_{0} \; .
                   \end{array}
   \right.
\EE 
The projection property of $J^{\bot}$ follows easily from that of $J$. }

The projector $J$ will now be used to establish the Lemma which 
describes the one-dimensional differential calculi associated with 
$\Gamma_{0}$, respectively $\hat{\Gamma}_{0}$.

\noindent {\bf Lemma 10: }{\sl Let ${\cal A}$ be a quantum group 
and $\Gamma_{0}$ , respectively $\hat{\Gamma}_{0}$, bimodules as defined 
in eqn.(33), respectively eqn.(86). Then we have

\noindent -(1) The one-dimensional bimodules $\Gamma_{0}$ and 
$\hat{\Gamma}_{0}$ are isomorphic.

\noindent -(2) The bimodule $\Gamma_{0}$ can be completed to a first 
order differential calculus over ${\cal A}$, with the differential 
$\delta$. This  differential calculus 
$(\Gamma_{0} , \delta)$ is isomorphic to the restriction 
$(\hat{\Gamma}_{0} , d|_{\hat{\Gamma}_{0}})$ of the Inner 
Differential Calculus $(\hat{\Gamma} , d)$ on its one-dimensional 
sub-bimodule $\hat{\Gamma}_{0}$. } 

\noindent {\bf Proof:} (1) The one-dimensional bimodule $\Gamma_{0}$ is 
defined as a free left-${\cal A}$-module generated by the element $X$. 
This element obeys the characteristic condition $X \wedge X = 0$ from 
eqn.(39). The right-multiplication of any $a \in{\cal A}$ is given by 
the restriction of the eqns.(35),(36) to the subspace $\Gamma_{0}$ 
\BE
Xa -aX = \partial|_{\Gamma_{0}} \equiv 0 \; .
\EE
This is an implication of the fact that the differential $\partial$ has 
no image in $\Gamma_{0}$ (${\rm Im}(\partial) \not\supset \Gamma_{0}$). 
The eqn.(94) can thus be replaced by the suggestive form
\BE
Xa - (\epsilon_{c} * a)X = 0
\EE
using the fact that ($\epsilon_{c} * a) = a$) for any $a \in {\cal A}$.

The bimodule $\hat{\Gamma}_{0}$ is the restriction of the bimodule 
$\hat{\Gamma}$, given by eqn.(82), to the one-dimensional 
sub-bimodule $\hat{\Gamma}_{0}$, which is generated by $\hat{X} = 
\sum_{a} \omega_{a}^{\;\; a}$, the canonical element of the bimodule 
 $\hat{\Gamma}$. This  sub-bimodule $\hat{\Gamma}_{0}$ has, up to the 
right-multiplication rule of elements $a \in {\cal A}$
\BE
\hat{X}a - (f_{a\;\;\;\;\; b}^{\;\; ab} * a)\hat{X} = 0 \; ,
\EE 
by construction, the same properties as the bimodule $\Gamma_{0}$, i.e. 
we have 
\BE
\hat{X} \wedge \hat{X} = 0 = X \wedge X \; .
\EE

The comparison of the eqns.(94) and (95) shows that one obtains an 
isomorphism of the left-pairs $(a , X)$ and $(a ,\hat{X} )$ by simply 
identifying the generating elements $X$ and $\hat{X}$. Concerning the 
associated right-pairs  $(X , b)$ and $(\hat{X} , a)$ the isomorphism 
is given by the relation
\BE
b = \sum_{a,b} (f_{a\;\;\;\;\; b}^{\;\; ab} * a) \; ,
\EE
that is, the isomorphism is the convolution product of any 
$a \in {\cal A}$ with the characteristic functionals given in eqn.(83), 
and the identification of the generating elements $X$ and $\hat{X}$. In 
this case the isomorphism is relative to the induced left-pairs 
$(b , X)$ and $((f_{a\;\;\;\;\; b}^{\;\; ab} * a) , \hat{X})$.

(2) The  right-multiplication rule described in eqn.(35) is due to the 
above given isomorphism of $\Gamma_{0}$ and $\hat{\Gamma}_{0}$ to be 
interpreted as the simplest choice of a right-multiplication rule. Thus 
we can define by an appropriate choice of a right-multiplication rule 
(we add simply a term $\delta a\sim X$) a differential $\delta$ on the 
bimodule $\Gamma_{0}$ 
\BE
0 \stackrel{{\rm def}}{=} Xa - (f_{\;\; 0}^{0} * a)X \; ,
\EE
where $f_{\;\; 0}^{0}$ generalizes the counit ${\epsilon}_{c}$ in 
eqn.(95). Thus we have analogously to eqn.(40)
\BEA
\delta a &=& Xa - aX  \nonumber \\
 & &\\
\delta \Theta  &=& \left[ X , \Theta \right]_{grad} |_{\Gamma_{0}} \; ,
\hspace{3cm}  \forall \; \Theta \, \in \Gamma_{0}^{\wedge}
\; ,  \nonumber 
\EEA
with the graded Lie-bracket
\begin{center}
$\left[ \;\;  , \; \right]_{grad} : \Gamma_{0} \times \Gamma_{0}^{\wedge} 
\longrightarrow \Gamma_{0}^{\wedge}$~, 
\end{center}
and if one constructs, as in theorem 7, a vector field $\chi_{X}$ which 
is dual to $X$, then one has 
\begin{center}
$\delta a = (\chi_{X} * a)X \; $~.
\end{center}
This makes $(\Gamma_{0} , \delta)$ a bicovariant differential calculus 
on ${\cal A}$.

The restriction of the Inner Differential Calculus $(\hat{\Gamma} , d)$ 
on its one-dimensional sub-bimodule $\hat{\Gamma}_{0}$, due to the 
projection operator $J$ given in definition 9, gives
\BE
J: \hat{\Gamma} \longrightarrow \hat{\Gamma}_{0} \hspace{3cm} 
J \circ d = d|_{\hat{\Gamma}_{0}} \; ,
\EE
which is an  one-dimensional differential calculus 
$(\hat{\Gamma}_{0}  , J \circ d)$ over ${\cal A}$. Due to the fact that 
the characteristic functionals $f_{\;\; 0}^{0}$ and 
$f_{a\;\;\;\;\; b}^{\;\; ab}$ involved in the eqns.(96),(99) are 
elements of the Hopf-algebra ${\cal A}'$ dual to the quantum group 
${\cal A}$, there must exist inverse functionals (in the sense of the 
antipode) of them.Thus we can define functionals 
$\tilde{f}_{\;\; 0}^{0}$ such that 
\begin{center}
$\tilde{f}_{\;\; 0}^{0}: f_{a\;\;\;\;\; b}^{\;\; ab} \longrightarrow 
f_{\;\; 0}^{0}$~.
\end{center}
Explicitly, this means that for all pairs $f_{\;\; 0}^{0} , 
f_{a\;\;\;\;\; b}^{\;\; ab} \in {\cal A}'$ there exists a
\BE
\tilde{f}_{\;\; 0}^{0} \stackrel{{\rm def}}{=} 
f_{\;\; 0}^{0} \, \kappa_{d}(f_{a\;\;\;\;\; b}^{\;\; ab}) \; ,
\EE
which obeys the relation 
\begin{center}
$m_{d}(\tilde{f}_{\;\; 0}^{0} \otimes f_{a\;\;\;\;\; b}^{\;\; ab}) = 
\tilde{f}_{\;\; 0}^{0} f_{a\;\;\;\;\; b}^{\;\; ab} = 
f_{\;\; 0}^{0} \kappa_{d}(f_{a\;\;\;\;\; b}^{\;\; ab}) 
f_{a\;\;\;\;\; b}^{\;\; ab} = f_{\;\; 0}^{0}$~. 
\end{center}
The converse mapping from $f_{\;\; 0}^{0}  \, \rightarrow \,  
f_{a\;\;\;\;\; b}^{\;\; ab}$ is defined analogously. Observing that the 
the characteristic functionals $f_{\;\; 0}^{0}$ and 
$f_{a\;\;\;\;\; b}^{\;\; ab}$ respectively fix, by the eqn.(99) 
respectively eqn.(96), the respective differentials $\delta$ and $J 
\circ d$ of the bimodules $\Gamma_{0}$ and $\hat{\Gamma}_{0}$, we see 
due to eqn.(102), the possibility to rewrite the differentials $\delta$ 
and $J \circ d$ into one another.

The combination of the results of the first part of our proof, that is 
the isomorphism of the bimodules $\Gamma_{0}$ and $\hat{\Gamma}_{0}$, 
and the possibility to rewrite $J \circ d$ into $\delta$ and vice versa 
makes the one-dimensional differential calculi $(\hat{\Gamma}_{0}  , 
J \circ d)$ and $(\Gamma_{0} , \delta)$ over ${\cal A}$ isomorphic.
\hfill \rule{3mm}{3mm}
\vspace{3mm}

{\sl Remark:} A realization of the  characteristic functionals 
$f_{\;\; 0}^{0}$, given in eqn.(99) in the proof of our Lemma, is 
given for example by the functional $f_{z}$, which is defined in Thm.5.6 
of \cite{Wo3}. This functional has all the properties, which one has to 
have for the characteristic functionals $f_{ij}$, as setteled in the 
eqns.(15)-(18) and in theorem 3. That is, the funcional
 
\newpage 
\noindent has to obey
\begin{center}
$ f_{z}(I) = 1$   
\end{center}
and
\begin{center}
$ f_{z}(ab) = f_{z}(a)f_{z}(b) $~,
\end{center}
for any  $a,b \in {\cal A}$. This funcionals approach in the limit $ z 
\rightarrow 0$ the counit ${\epsilon}_{c}$, that is one recovers 
eqn.(95) in this limit. A further example of such functionals (in the 
case of $SL_{q}(N)$ quantum groups) is described in \cite{SS1}, in the 
remark 4 following theorem 2.2.

In order to map the Outer Differential Calculi $(\Gamma ,\partial)$ onto 
an Inner Differential Calculi $(\hat{\Gamma} , d)$ one has to extend the 
bimodule $\Gamma$ by the bimodule which is characterized by the 
generating element of the exterior differential $\partial$, i.e. by 
$\Gamma_{0}$ . Thus we obtain the extended  bimodule $\tilde{\Gamma}$, 
given in eqn.(33)
\begin{center}
$\tilde{\Gamma} = \Gamma \oplus \Gamma_{0}$~.
\end{center}
Furthermore we have to complete the differential $\partial$ in such a 
manner, that the image of this extended differential $\partial'$, 
${\rm Im}(\partial')$, with ${\rm Im}(\partial) \subset 
{\rm Im}(\partial')$ and    ${\rm Im}(\partial')/{\rm Im}(\partial) 
\sim X$ contains a term proportional to $X$. This is obtained due to a 
modification of the right-multiplication rule eqn.(35), in the following 
way
\BEA
\tilde{\xi}a &\stackrel{{\rm def}}{=}& c (f_{\;\; 0}^{0} * a)X + 
c(\partial a + \xi a)  \nonumber \\
\hspace{3cm} &=&  c \left[ (f_{\;\; 0}^{0} - \epsilon_{c}) * a 
\right]X  + c(\epsilon_{c} * a)X +c(\partial a + \xi a) \; ,\hspace{2cm}
\EEA
where $\tilde{\xi}$ and $\xi$ are defined as in eqn.(35). The functional 
$f_{\;\; 0}^{0} \in {\cal A}_{d}$, is as in eqn.(99) the characteristic 
functional of the $({\cal A} - \Gamma_{0})$-commutation, ruling the 
commutation of the elements $a \in {\cal A}$ with the generator $X$ of 
the bimodule $\Gamma_{0}$.

Observe that eqn.(103) becomes identical to eqn.(35) in the case 
$f_{\;\; 0}^{0} = \epsilon_{c}$. Here seems to be the right place to 
remark that the counit $\epsilon_{c}$, at least on the formal level, has 
also the properties of the characteristic functionals $f_{ij}$, as 
described in the conditions eqns.(15)-(18). But the vector field 
associated to such a differential calculus is degenerate (i.e. the 
vector field is identically 0). Thus, in a certain sense, we can think 
about Outer Differential Calculi as Inner Differential Calculi, which 
are degenerate in the sector generating the differential.

Following now the argumentation, which leads from eqn.(35) to eqn.(36), 
we obtain, applying it to eqn.(103)
\BE
Xa - (\epsilon_{c} * a)X = \partial a + 
\underbrace{\left[ (f_{\;\; 0}^{0} - \epsilon_{c}) * a 
\right]X}_{\delta a} = \partial a + \delta a\; ,
\EE
where we have set $ \delta a = (f_{\;\; 0}^{0} - \epsilon_{c}) * aX = 
(\chi_{X} * a)X$. Here $\chi_{X}$ is the vectorfield dual to the 
generating element $X$ of  $\Gamma_{0}$, analogously to the constructions 
subjet of theorem 6. (Compare the definition of $\chi_{X}$ with 
eqn.(78) in the case $\chi_{X} \equiv \chi_{\;\; c}^{c}$.)

As in the eqns.(40)-(42) we obtain for any $\tilde{\Theta} \in 
\tilde{\Gamma}^{\wedge}$
\BE
(\partial + \delta)\tilde{\Theta} = \left[ X \, , \, \tilde{\Theta}
\right]_{grad} \; ,
\EE
where $ \left[ X \, , \; \right]_{grad}$ now maps onto the whole of  
$\tilde{\Gamma}^{\wedge}$
\BE
\left[ X \, , \, \, \right]_{grad} : \tilde{\Gamma}^{\wedge n}
\longrightarrow  \tilde{\Gamma}^{\wedge n + 1} \, .
\EE
This is nothing but the encoding of our right-multiplication rule (103) 
in the graded bracket.

Due to this construction we get the extended Outer Differential 
Calculi\footnote{ This differential calculi are Inner Differential 
Calculi in the sense of our classification.} 
$(\tilde{\Gamma} \, , \, \partial + \delta) = (\Gamma \, , \, \partial) 
\oplus (\Gamma_{0} \, , \, \delta)$ over ${\cal A}$, this is the 
completion of the Outer Differential Calculi $(\Gamma \, , \, \partial)$ 
by some complementary one-dimensional differential calculus 
$(\Gamma_{0} \, , \, \delta)$ as given in Lemma 10.

This differential calculi $(\tilde{\Gamma} \, , \, \partial + \delta)$ 
are isomorphic to Inner Differential Calculi $(\hat{\Gamma} , d)$, 
which follows from the isomorphy of $(\Gamma_{0} \, , \, \delta)$ and 
the one-dimensional sub- differential calculi $(\hat{\Gamma}_{0} , 
J \circ d)$ shown in Lemma 10, and Proposition 2. Because, due to 
Proposition 2, the part of the Inner  Differential Calculi 
$(\hat{\Gamma} , d)$ which is complementary to $(\hat{\Gamma}_{0} , 
J \circ d)$ (that is $(J^{\bot} \circ \hat{\Gamma} , J^{\bot} 
\circ d)$), which is defined by the projection operator $J^{\bot}$ 
given in eqn.(93), is isomorphic to some Outer Differential Calculi 
$(\Gamma \, , \, \partial)$. Thus we have shown the isomorphism of the 
differential calculi $(\tilde{\Gamma} \, , \, \partial + \delta)$ and a 
certain Inner  Differential Calculi $(\hat{\Gamma} , d)$.

This construction allows us to give a family of natural, injective maps
\BE
\Delta^{Out}_{In}: \{(\Gamma ,\partial)\} \longrightarrow 
\{(\hat{\Gamma} , d)\} \; ,
\EE
which maps the set of isomorphy-classes $\{(\Gamma ,\partial)\}$ of 
Outer Differential Calculi into the set of isomorphism-classes 
$\{(\hat{\Gamma} , d)\}$ of Inner  Differential Calculi. These maps are 
parametrized by the characteristic functionals $f_{\;\; 0}^{0}$, or 
equivalently by the choice of $\delta$, as was shown in the eqns.(103),
(104).

The map of an arbitrarily chosen, but then fixed, Outer Differential 
Calculus $(\Gamma \, , \, \partial)$ is given by
\BE
\Psi_{\Gamma  , \Delta^{Out}_{In}(\Gamma )}:  (\Gamma \, , \, \partial) 
\longrightarrow \Delta^{Out}_{In}(\Gamma \, , \, \partial) \; .
\EE
This map is given by 
\BE
\Psi_{\Gamma  , \Delta^{Out}_{In}(\Gamma )}:  \Gamma \longrightarrow 
\Gamma \oplus \Gamma_{0} 
\EE
and 
\BE
\Psi_{\Gamma  , \Delta^{Out}_{In}(\Gamma )}:  \partial \longrightarrow 
\partial  \oplus \delta \; ,
\EE
as in our previous construction.

The map from an Inner Differential Calculus $(\hat{\Gamma} , d)$ onto an 
Outer Differential Calculus $(\Gamma \, , \, \partial)$ is obtained by 
the restriction of $(\hat{\Gamma} , d)$ on the part which is 
complementary to $(\hat{\Gamma}_{0} , J \circ d)$. This is done with the 
help of the projection operators given in def. 9. Applying $J^{\bot}$ 
from eqn.(93) on $\hat{\Gamma}$ we obtain 
\BE
J^{\bot} \hat{\Gamma} = \hat{\Gamma}/\hat{\Gamma}_{0} \; ,
\EE
which is the sub-bimodule of $\hat{\Gamma}$ spanned by that 
basis-elements of $\hat{\Gamma}$, which are the orthocomplement of the 
generator $\hat{X}$. Thus we have the situation which we have met in the 
case of Outer Differential Calculi.

Now we split the differential $d$ into two pieces, such that the first 
part will map onto the sub-bimodule $\hat{\Gamma}/\hat{\Gamma}_{0}$, and 
the second part maps onto the complementary  sub-bimodule 
$\hat{\Gamma}_{0}$. This splitting is given by the projectors $J$ and 
$J^{\bot}$ defined in eqns.(92),(93), which complete each other to the 
identity-map on $\hat{\Gamma}$ ($(id_{\hat{\Gamma}} = J + J^{\bot})$) 
\BE
d = (J + J^{\bot}) \circ d \; .
\EE
With the properties given in Def. 9 we have the maps
\BEA
J \circ d & \hat{\Gamma}^{\wedge n} \longrightarrow 
  \hat{\Gamma}^{\wedge n} \wedge \hat{\Gamma}_{0} & 
  \subset \hat{\Gamma}^{\wedge n +1}  \\
  J^{\bot} \circ d & \hat{\Gamma}^{\wedge n} \longrightarrow 
  \hat{\Gamma}^{\wedge n} \wedge (\hat{\Gamma}/ \hat{\Gamma}_{0}) &
  \subset \hat{\Gamma}^{\wedge n +1} \; .
\EEA
Then we can define the exterior differential 
\BE
(J^{\bot} \circ d)\; \hat{\Theta} = J^{\bot} \circ \frac{1}{\lambda} 
\left[ \hat{X} \, , \, \hat{\Theta} \right]_{grad} = 
\frac{1}{\lambda} \left[ \hat{X} \, , \, 
\hat{\Theta}\right]_{grad}|_{\hat{\Gamma}^{\wedge}/ \hat{\Gamma}_{0}}\; ,
\EE
on $\hat{\Gamma}$, where $\hat{\Theta} \in \hat{\Gamma}^{\wedge}$ and 
$\left[ \; , \, \right]_{grad}$ corresponds to the graded Lie-bracket 
given in eqn.(87), while $\hat{\Gamma}^{\wedge}/ \hat{\Gamma}_{0}$ is 
the part of the exterior algebra over $\hat{\Gamma}$ which contains no 
$\hat{\Gamma}_{0}$-factors. This implies readily that eqn.(115) has the 
same structure like eqn.(40). Thus it makes the differential calculus 
$(J^{\bot} \circ \hat{\Gamma} , J^{\bot} \circ d)$ an Outer Differential 
Calculus. Remembering now the universality, stated in Prop. 2, there 
must exist an  Outer Differential Calculus $(\Gamma , \partial)$ which 
is isomorphic to $(J^{\bot} \circ \hat{\Gamma} , J^{\bot} \circ d)$.

The considerations from above give us a natural, surjective map 
\BE
\Delta^{In}_{Out}: \{(\hat{\Gamma} , d)\} \longrightarrow 
\{(\Gamma ,\partial)\} \; ,
\EE
which maps the set of isomorphism-classes $\{(\hat{\Gamma} , d)\}$ of 
Inner  Differential Calculi onto the set of isomorphism-classes 
$\{(\Gamma ,\partial)\}$ of Outer Differential Calculi. 

For an arbitrarily chosen, but then fixed, Inner Differential Calculus 
$(\hat{\Gamma} , d)$ the map
\BE
\Phi_{\hat{\Gamma} , \Delta^{In}_{Out}(\hat{\Gamma})}: 
(\hat{\Gamma} , d) \longrightarrow (\Gamma , \partial)
\EE 
associates to it an Outer Differential Calculus $(\Gamma , \partial)$ by 
\BE
\Phi_{\hat{\Gamma} , \Delta^{In}_{Out}(\hat{\Gamma})}: \hat{\Gamma} 
\longrightarrow J^{\bot} \circ \hat{\Gamma} = 
\hat{\Gamma}/ \hat{\Gamma}_{0} \equiv \Gamma
\EE
and 
\BE
\Phi_{\hat{\Gamma} , \Delta^{In}_{Out}(\hat{\Gamma})}: d \longrightarrow 
J^{\bot} \circ d \equiv 
\partial \; ,
\EE
given due to the projector, like in our construction above.

Now we are prepared to state the reconstruction theorem.

\noindent {\bf Theorem 11:} {\sl Let ${\cal A}$ be a quantum group, 
which allows for the construction of Inner- as well as Outer 
Differential Calculi. Then we have

\noindent -(1) There exists a natural, injective family of maps
\begin{center}
$\Delta^{Out}_{In}: \{(\Gamma ,\partial)\} \longrightarrow 
\{(\hat{\Gamma} , d)\}$~,
\end{center}
which maps the set of isomorphy-classes $\{(\Gamma ,\partial)\}$ of 
Outer Differential Calculi into the set of isomorphy-classes 
$\{(\hat{\Gamma} , d)\}$ of Inner  Differential Calculi. 

For an arbitrarily chosen, but then fixed, Outer Differential Calculus 
$(\Gamma \, , \, \partial)$, the map
\begin{center}
$ \Psi_{\Gamma  , \Delta^{Out}_{In}(\Gamma )}:  
(\Gamma \, , \, \partial) \longrightarrow 
(\Gamma \oplus \Gamma_{0} \, , \, \partial \oplus \delta ) \; .$
\end{center}
is an ${\cal A}$-bimodule-monomorphism of  $\Gamma$, which is 
parametrized by $f_{\;\; 0}^{0}$. This parametrization fixes also the 
differential $\delta$. 

\noindent -(2) There exists a natural, surjective map 
\begin{center}
$\Delta^{In}_{Out}: \{(\hat{\Gamma} , d)\} \longrightarrow 
\{(\Gamma ,\partial)\}$~,
\end{center} 
which maps the set of isomorphy-classes $\{(\hat{\Gamma} , d)\}$ of 
Inner  Differential Calculi onto the set of isomorphy-classes 
$\{(\Gamma ,\partial)\}$ of Outer Differential Calculi. 

For an arbitrarily chosen, but then fixed, Inner Differential Calculus 
$(\hat{\Gamma} , d)$, the map
\begin{center}
$\Phi_{\hat{\Gamma} , \Delta^{In}_{Out}(\hat{\Gamma})}: 
(\hat{\Gamma} , d) \longrightarrow 
(\Gamma , \partial)$ 
\end{center} 
is an ${\cal A}$-bimodule-epimorphism of $\hat{\Gamma}$. The 
differential $\partial$ is fixed by the given differential $d$.

\noindent -(3) On the set of isomorphy-classes  
$\{(\Gamma ,\partial)\}$ of Outer Differential Calculi we have the 
relation
\begin{center}
$id_{\{(\Gamma ,\partial)\}} = \Delta^{In}_{Out} \circ 
\Delta^{Out}_{In}$~. 
\end{center}  }

\noindent {\bf Proof:} (1) The family $\Delta^{Out}_{In}$ is given in 
eqn.(107), due to the construction of an extended bimodule 
$\tilde{\Gamma}$ and an associated extension of the differential given in 
eqns.(103),(104) by a modified right-multiplication rule.

The map $\Psi_{\Gamma  , \Delta^{Out}_{In}}$ given for an arbitrarily 
fixed Outer Differential Calculus $(\Gamma \, , \, \partial)$ is  
parametrized by the characteristic functional $f_{\;\; 0}^{0}$, ruling 
the ($\Gamma_{0} - {\cal A}$)-commutation. It is fixed in the 
eqns.(108)-(110). The fact that 
\BE
ker(\Psi_{\Gamma  , \Delta^{Out}_{In}(\Gamma )}) = \{\emptyset \}
\EE
is due to the property of mapping the Outer Differential Calculus 
$(\Gamma \, , \, \partial)$ identically onto itself.

Now let $\Gamma_{i}$ ($ i \in \cal I$ , $\cal I$ an indexset) denote the 
representatives of the isomorphy-classes $\{(\Gamma ,\partial)\}$, 
then we have by eqn.(109)
\BE
\Psi_{\Gamma  , \Delta^{Out}_{In}(\Gamma )}( \Gamma_{i}) = 
\Gamma_{i} \oplus \Gamma_{0} \equiv \hat{\Gamma}_{i} \; . 
\EE
Thus, because the $\Gamma_{i}$ are by assumption from different 
isomorphy-classes, the $\hat{\Gamma}_{i} $ are in different 
isomorphy-classes too, and have all a unique $\Gamma_{i}$ as inverse 
image. Thus the maps $\Psi_{\Gamma  , \Delta^{Out}_{In}(\Gamma )}$ 
are ${\cal A}$-bimodule-monomorphisms of the bimodules $\Gamma$. The 
image ofthe differential $\partial$ by $\Psi_{\Gamma  , 
\Delta^{Out}_{In}(\Gamma )}$ is uniquely fixed by the functionals 
$f_{\;\; 0}^{0}$. This finishes the proof of (1).

\noindent (2) The map $\Delta^{In}_{Out}$ is given in eqn.(116) and is 
fixed due to the projector $J^{\bot}$ given in eqn.(93), which projects 
$\hat{\Gamma}$, as well as $d$ to the corresponding objects ($\Gamma$ 
and $\partial$). The surjectivity is seen due to the fact that, with the 
map  $\Delta^{Out}_{In}$ we can allways get an Inner Differential 
Calculus, which is mapped to a prescribed Outer Differential Calculus 
$(\Gamma \, , \, \partial)$ by  $\Delta^{In}_{Out}$.

The map $\Phi_{\hat{\Gamma} , \Delta^{In}_{Out}(\hat{\Gamma})}$ applied 
to any Inner Differential Calculus $(\hat{\Gamma} , d)$ given in the 
eqns. (117)-(119) has a nontrivial kernel 
\BE
ker(\Phi_{\hat{\Gamma} , \Delta^{In}_{Out}(\hat{\Gamma})}) = 
(\hat{\Gamma}_{0} , \delta) \; .
\EE
Its ${\cal A}$-bimodule-epimorphism property follows from the fact that, 
having any Outer Differential Calculus $(\Gamma \, , \, \partial)$, one 
can with the application of the map $\Psi_{\Gamma  , \Delta^{Out}_{In}
(\Gamma )}$ associate an Inner Differential Calculus 
$(\tilde{\Gamma} \, , \, \partial + \delta)$. The application of the map 
$\Phi_{\hat{\Gamma} , \Delta^{In}_{Out}(\hat{\Gamma})}$ on this Inner 
Differential Calculus gives us back the calculus we have started with. 
Thus it is possible to obtain any Outer Differential Calculi, 
using that $\Phi_{\hat{\Gamma} , \Delta^{In}_{Out}(\hat{\Gamma})}$ is an  
${\cal A}$-bimodule-epimorphism. Furthermore we have
\BE
\Phi_{\hat{\Gamma} , \Delta^{In}_{Out}(\hat{\Gamma})} \circ 
\Psi_{\Gamma  , \Delta^{Out}_{In}(\Gamma )}(\partial) = \partial \; ,
\EE
and we have finished the proof of (2). 

\noindent (3) The set of isomorphy-classes $\{(\Gamma ,\partial)\}$ 
of Outer Differential Calculi is mapped by $\Delta^{Out}_{In}$ into the 
set of isomorphy-classes  $\{(\hat{\Gamma} , d)\}$ of Inner 
Differential Calculi as was shown in (1). The application of the map 
$\Delta^{In}_{Out}$ given in (2) on an  Inner Differential Calculus 
obtained as above gives obviously, due to the projection property, back 
the original calculus we started with. Thus we have verified the relation
\BE
id_{\{(\Gamma ,\partial)\}} = \Delta^{In}_{Out} \circ 
\Delta^{Out}_{In} \; ,
\EE 
which finishes the proof. 
\hfill \rule{3mm}{3mm}
\vspace{3mm}

{\sl Remark:} All our considerations which lead to Theorem 11 can be 
applied to left- as well as right-invariant differential calculi.
\vspace{15mm}

\noindent {\Large {\bf 3. A Differential Bicomplex on Quantum Groups}}
\vspace{6mm}

The eqn.(87) taking into account the results obtained in Lemma 10 and 
Theorem 11 allows for a splitting of the Cartan condition 
\BE
0 = d^{2} = (\partial + \delta)^{2} = \partial^{2} + \partial\delta 
+ \delta\partial + \delta^{2} \; .
\EE
Remembering the relation  
\BE
\partial^{2} = 0 \; ,
\EE
that is eqn.(43), and using a combination of the eqns.(97) and (100) to 
obtain the identity 
\BE
\delta^{2} = 0 \; ,
\EE
we have at hand three independent Cartan conditions. Thus one obtains, 
by the insertion of the eqns.(126) and (127) in eqn.(125) the relation
\BE
\partial\delta + \delta\partial = 0 \; .
\EE 

The splitting of the exterior differential $d = (\partial + \delta)$ and 
the corresponding decomposition of the Cartan condition eqn.(125) allows 
us to rewrite the q-de Rham complex 
\BE
0 \longrightarrow {\cal A} \stackrel{d}{\longrightarrow} 
\hat{\Gamma}^{\wedge 1} \stackrel{d}{\longrightarrow} 
\hat{\Gamma}^{\wedge 2} \stackrel{d}{\longrightarrow} 
\cdots \stackrel{d}{\longrightarrow} 
\hat{\Gamma}^{\wedge ({\rm dim}\Gamma + 1)}\stackrel{d}{\longrightarrow} 
0 \; ,
\EE
of the Inner Differential Calculus $(\hat{\Gamma} , d)$ in a 
differential bicomplex $(\tilde{\Gamma} , \partial , \delta)$ 
\BE
\begin{array}{lclclclclcl}
\tilde{\Gamma}^{0,0}&\stackrel{\partial}{\longrightarrow}
&\tilde{\Gamma}^{0,1}&\stackrel{\partial}{\longrightarrow}
&\tilde{\Gamma}^{0,2}&\stackrel{\partial}{\longrightarrow}
&\cdots&\stackrel{\partial}{\longrightarrow}
&\tilde{\Gamma}^{0,{\rm dim}\Gamma}&\stackrel{\partial}{\longrightarrow}
&0 \\
\bigg\downarrow {\mbox{\footnotesize $\delta$}}& &
\bigg\downarrow {\mbox{\footnotesize $\delta$}}& 
&\bigg\downarrow {\mbox{\footnotesize $\delta$}}& & & 
&\bigg\downarrow {\mbox{\footnotesize $\delta$}}& & \\
\tilde{\Gamma}^{1,0}&\stackrel{\partial}{\longrightarrow}
&\tilde{\Gamma}^{1,1}&\stackrel{\partial}{\longrightarrow}
&\tilde{\Gamma}^{1,2}&\stackrel{\partial}{\longrightarrow}
&\cdots&\stackrel{\partial}{\longrightarrow}
&\tilde{\Gamma}^{1,{\rm dim}\Gamma}&\stackrel{\partial}{\longrightarrow}
&0 \\
\bigg\downarrow {\mbox{\footnotesize $\delta$}}& 
&\bigg\downarrow {\mbox{\footnotesize $\delta$}}& 
&\bigg\downarrow {\mbox{\footnotesize $\delta$}}& & & 
&\bigg\downarrow {\mbox{\footnotesize $\delta$}}& & \\
0& &0& &0& & & &0& & \; .
\end{array}
\EE
where we have used the notation $\tilde{\Gamma}^{0,0} = {\cal A} , 
\tilde{\Gamma}^{0,1} = \Gamma_{0}$ and $\tilde{\Gamma}^{1,0} = \Gamma$ 
from eqn.(33). The spaces $\hat{\Gamma}^{\wedge q +1}$ of the q-de Rham 
complex are decomposed into two complementary spaces 
$\tilde{\Gamma}^{0,q + 1}$ and $\tilde{\Gamma}^{1,q}$. This reflects 
the splitting of the spaces of $q + 1$-forms 
$\hat{\Gamma}^{\wedge q +1} = \tilde{\Gamma}^{0,q + 1} 
\oplus \tilde{\Gamma}^{1,q}$ into parts which are zero-forms with respect 
to $\Gamma_{0}$ (in the first case), respectively one-forms with respect 
to $\Gamma_{0}$ (in the second case).
\vspace{4mm}

{\small \it Acknowledgement: I would like to thank Prof. J.Wess for the 
invitation to Munich and the Max-Planck-Institut f{\"u}r Physik in 
Munich for financial support. Dr. M.Schlichenmaier is acknowledged for 
helpful discussions.}

\newpage 
\end{document}